\documentclass[a4paper,openany, 12pt]{article}

\usepackage[small]{titlesec}

\usepackage{booktabs}
\usepackage[T1]{fontenc}
\usepackage[utf8]{inputenc}
\usepackage[english]{babel}
\usepackage[a4paper,top=4cm,bottom=3cm,left=2.5cm,right=2.5cm]{geometry}
\usepackage{amsfonts}
\usepackage{mathtools}
\usepackage{amsmath}

\usepackage{lmodern}
\usepackage{empheq}
\usepackage{amsthm}
\usepackage{amssymb}
\usepackage{booktabs}
\usepackage{caption}
\usepackage{siunitx}
\usepackage{float}
\usepackage{amssymb}
\usepackage{graphicx}
\usepackage{stmaryrd}
\usepackage{color}
\usepackage{pifont}%
\usepackage{enumerate}
\usepackage{enumitem}
\setlist[description]{leftmargin=\parindent,labelindent=\parindent}
\usepackage{listings} 
\usepackage{amsmath}
\usepackage{amscd}
\usepackage{braket}
\usepackage{physics}
\usepackage{tikz}
\usepackage{calc}
\usepackage[hidelinks]{hyperref}

\usepackage{todonotes}

\newcommand{\dom}{\mathrm{dom}}

\newcommand{\R}{\mathbb{R}}
\newcommand{\sph}{\mathbb{S}}
\newcommand{\K}{\mathcal{K}}

\newcommand{\haus}{\mathcal{H}}
\newcommand{\Coco}{\Conv_{\mathrm{cd}}(\mathbb{R}^n)}

\newtheorem{theorem}{Theorem}[section]
\newtheorem{corollary}[theorem]{Corollary}

\newtheorem{lemma}[theorem]{Lemma}

\newcommand*\sq{\mathbin{\vcenter{\hbox{\rule{.3ex}{.3ex}}}}}
\NewDocumentCommand{\sff}{}{\mathrm{I\!I}}

\newcommand{\Conv}{\mathrm{Conv}}
\newcommand{\vol}{\mathrm{vol}}

\newcommand{\supp}{\mathrm{supp}\ }

\newcommand{\GW}{\mathrm{GW}}

\newcommand{\Cosc}{\mathrm{Conv}_{\mathrm{sc}}}

\newcommand{\epi}{\mathrm{epi}}

\newcommand\blfootnote[1]{%
  \begingroup
  \renewcommand\thefootnote{}\footnote{#1}%
  \addtocounter{footnote}{-1}%
  \endgroup
}
\title{From valuations on convex bodies to convex functions}
\author{Jonas Knoerr and Jacopo Ulivelli}
\date{}

\begin{document}
\maketitle

\begin{abstract}
    A geometric framework relating valuations on convex bodies to valuations on convex functions is introduced. It is shown that a classical result by McMullen can be used to obtain a characterization of continuous, epi-translation invariant, and $n$-epi-homogeneous valuations on convex functions, which was previously established by Colesanti, Ludwig, and Mussnig. Following an approach by Goodey and Weil, a new characterization of $1$-epi-homogeneous valuations is obtained.
    
\end{abstract}

\section{Introduction}

For a family $\mathcal{S}$ of subsets of $\R^n$, a functional $Y: \mathcal{S} \to \R$ is a real-valued \textit{valuation} if for every $A, B \in \mathcal{S}$, 
\begin{equation}\label{valuation_property}
    Y(A\cap B)+Y(A \cup B)=Y(A)+Y(B)
\end{equation}
whenever $A\cup B, A \cap B \in \mathcal{S}$. This concept goes back to Dehn's solution of Hilbert's third problem and has since then played a central role in convex and discrete geometry (see {\cite[Chapter 6]{SchneiderConvexBodiesBrunn2013}} for a comprehensive overview of the subject). Valuations on convex bodies, that is, valuations on the space $\mathcal{K}^n$ of all non-empty, convex, and compact subsets of $\R^n$, have been the focus of intense research ever since, with many break-through results in recent years, see \cite{AleskerContinuousrotationinvariant1999,AleskerDescriptiontranslationinvariant2001,BernigFuHermitianintegralgeometry2011,BernigEtAlIntegralgeometrycomplex2014,HaberlMinkowskivaluationsintertwining2012,HaberlParapatitscentroaffineHadwiger2014,HaberlParapatitsMomentsvaluations2016,LudwigReitznerclassification$SLn$invariant2010}. 

Following these advances in modern valuation theory, the notion of valuations has been extended to families of functions. Here, a functional $Z: \mathcal{F} \to \R$ defined on a family $\mathcal{F}$ of extended real-valued functions is called a valuation if \[Z(u \wedge v)+Z(u \vee v)=Z(u)+Z(v)\] for all $u,v \in \mathcal{F}$ such that the pointwise minimum $u \wedge v$ and maximum $u \vee v$ belong to $\mathcal{F}$. If $\mathcal{F}$ denotes the family of indicator functions of convex bodies or the family of support functions $h_K$, given for $K\in\K^n$ by
\begin{align*}
    h_K(y)=\sup_{x\in K}x\cdot y, \quad y\in\R^n,
\end{align*} this recovers the classical notion of valuations on $\K^n$. In this sense, valuations on functions generalize valuations on sets. There is, however, also the following more geometric interpretation: Assume that the functions are defined on $\R^n$ and consider for $u\in\mathcal{F}$ its epi-graph 
\begin{align*}
	\epi(u)\coloneqq \{(x,t)\in\R^n\times\R: u(x)\le t \}\subset \R^{n+1}.
\end{align*}
It is easy to see that $\epi(u\wedge v)=\epi(u)\cup \epi(v)$, $\epi(u\vee v)=\epi(u)\cap\epi(v)$ for $u,v\in \mathcal{F}$. Thus, valuations on functions correspond to valuations on epi-graphs.\\

Due to their intimate relation with convex bodies, valuations on various spaces related to convexity have been considered, including log-concave and quasi-concave functions \cite{ColesantiLombardiValuationsspacequasi2017,ColesantiEtAlTranslationinvariantvaluations2018,MussnigValuationslogconcave2021}. However, there are also many results concerning other classical function spaces, for example, $L_p$ and Orlicz spaces \cite{KoneValuationsOrliczspaces2014,LiMaLaplacetransformsvaluations2017,TsangLp2010,TsangMinkowski2012}, continuous and Lipschitz \linebreak functions \cite{ColesantiEtAlclassinvariantvaluations2020,ColesantiEtAlContinuousvaluationsspace2021}, Sobolev spaces and spaces of functions of bounded variation \cite{LudwigFisherinformationmatrix2011,LudwigValuationsfunctionspaces2011,LudwigValuationsSobolevSpaces2012,MaRealvaluedvaluations2016,WangSemivaluations$BVRn$2014}, as well as definable functions \cite{BaryshnikovEtAlHadwigersTheoremdefinable2013}, or general Banach lattices \cite{TradaceteVillanuevaValuationsBanachlattices2020}.\\

Valuations on convex functions on $\R^n$ have been one of the most active areas of research \cite{AleskerValuationsconvexfunctions2019,CavallinaColesantiMonotonevaluationsspace2015,ColesantiEtAlMinkowskivaluationsconvex2017,ColesantiEtAlValuationsConvexFunctions2017,ColesantiEtAlHadwigertheoremconvex2020,ColesantiEtAlHessianValuations2020,ColesantiEtAlhomogeneousdecompositiontheorem2020,ColesantiEtAlHadwigertheoremconvex2021,ColesantiEtAlHadwigertheoremconvex2022,ColesantiEtAlHadwigertheoremconvex2022a,KnoerrSmoothvaluationsconvex,Knoerrsupportduallyepi2021,KnoerrSingularvaluationsHadwiger2022}. The primary focus is on subspaces of the space $\Conv(\R^n)$ of all convex functions $u:\R^n\rightarrow(-\infty,+\infty]$ that are lower semi-continuous and proper, that is, not identically $+\infty$. We will equip these spaces with the topology induced by epi-convergence (see Section \ref{section:ConvexFunctions} for details). In general, the family $\Conv(\R^n)$ seems to be too large to allow for a large class of interesting valuations; indeed, under seemingly reasonable invariance assumptions, there do not exist non-constant valuations on this space, as observed, for example, in {\cite[Section 9]{ColesantiEtAlhomogeneousdecompositiontheorem2020}}.

To obtain a sufficiently interesting set of valuations, one has to focus on smaller spaces of convex functions. One of the natural choices is the space of \emph{super coercive} convex functions used in \cite{ColesantiEtAlValuationsConvexFunctions2017,ColesantiEtAlHessianValuations2020,ColesantiEtAlhomogeneousdecompositiontheorem2020} \[ \Cosc(\R^n)\coloneqq \{u \in \Conv(\R^n): \lim_{|x|\to +\infty} u(x)/|x|=+\infty \},\] together with the space of \emph{finite} convex functions \[\Conv(\R^n,\R)\coloneqq \{v \in \Conv(\R^n): v(x)<+\infty \text{ for all }x\in\R^n \}.\]
These two families are in natural bijection using the \textit{Fenchel-Legendre transform} $u\mapsto u^*$, where \[u^{*}(y)=\sup_{x \in \R^n}  x\cdot y-u(x) =\sup_{(x,t)\in \epi(u)}  x\cdot y -t, \quad y\in\R^n. \]
It will be advantageous to think of elements of $\Cosc(\R^n)$ as unbounded convex sets in $\R^{n+1}$ via their epi-graph, whereas we will consider their Fenchel-Legendre transforms (that is, elements of $\Conv(\R^n,\R)$) as the associated support functions.  This leads to the following notion, introduced by Colesanti, Ludwig, and Mussnig in \cite{ColesantiEtAlhomogeneousdecompositiontheorem2020}: A valuation $Z:\Cosc(\R^n)\rightarrow\R$ is called \emph{epi-translation invariant} if
\begin{align*}
	Z(u(\cdot-x)+c)=Z(u)\quad\text{for all }u\in\Cosc(\R^n), x\in\R^n,c\in\R,
\end{align*}
that is, $Z$ is invariant with respect to translations of the epi-graph of $u$ in $\R^{n+1}$. Let us similarly define \textit{epi-multiplication} by
\begin{align*}
	t\sq u (x)\coloneqq t u\left(\frac{x}{t}\right)\quad \text{for }u\in\Cosc(\R^n), x\in\R^n, t>0,
\end{align*}
that is, the epi-graph of $t\sq u$ is obtained by rescaling $\epi(u)$ by $t>0$. We call a valuation $Z:\Cosc(\R^n)\rightarrow\R$ \textit{epi-homogeneous of degree $i$} if
\begin{align*}
	Z(t\sq u)=t^iZ(u)\quad\text{for all } u\in \Cosc(\R^n), t>0.
\end{align*}

Obviously, continuous epi-translation invariant valuations on $\Cosc(\R^n)$ are closely related to continuous and translation invariant valuations on convex bodies. As shown by Colesanti, Ludwig, and Mussnig, this relation can be exploited to establish a homogeneous decomposition theorem for this type of valuation, similar to the corresponding decomposition for translation invariant valuations on convex bodies obtained by McMullen \cite{McMullenValuationsEulerType1977}.
\begin{theorem}[Colesanti, Ludwig, and Mussnig \cite{ColesantiEtAlhomogeneousdecompositiontheorem2020}]\label{theorem:HomogeneousDecomposition}
If $Z: \Cosc(\R^n)\to\R$ is a continuous and epi-translation invariant valuation, then there are continuous and epi-translation invariant valuations $Z_0,\dots, Z_n: \Cosc(\R^n) \to \R$ such that $Z_i$ is epi-homogeneous of degree $i$ and $Z=Z_0+\dots +Z_n$.  
\end{theorem}
They were furthermore able to provide a complete characterization of valuations of maximal degree of homogeneity. Let $C_c(\R^n)$ denote the family of continuous functions on $\R^n$ with compact support. For a function $u:\R^n\to (-\infty,+\infty]$ we also define its \textit{domain} \[\dom(u)\coloneqq \{x\in \R^n: u(x)<+\infty \}.\]

\begin{theorem} [Colesanti, Ludwig, and Mussnig \cite{ColesantiEtAlhomogeneousdecompositiontheorem2020}]\label{theorem:CharMaxDegree} 
A functional $Z: \Cosc(\mathbb{R}^n) \to \mathbb{R}$ is a continuous and epi-translation invariant valuation
that is epi-homogeneous of degree $n$, if and only if there exists $\zeta \in C_c(\mathbb{R}^n)$ such that 
\begin{equation*}
    Z(u)=\int_{\dom(u)} \zeta(\nabla u(x))\, \,dx
\end{equation*}
for every $u \in \Cosc(\mathbb{R}^n)$.
\end{theorem}

In \cite{Knoerrsupportduallyepi2021}, a version of Theorem \ref{theorem:HomogeneousDecomposition} was obtained by interpreting valuations on convex functions as valuations on convex bodies in $\R^{n+1}$ using their support function via the map $$K \mapsto h_K(\cdot,-1).$$ Here, we follow a more geometric approach by considering the epigraphs of functions in $\Cosc(\R^n)$ as the "bodies" we are interested in. In order to make this precise, we consider the space 
\begin{align*}
    \Coco\coloneqq \{u\in \Cosc(\R^n): u \text{ has compact domain} \},
\end{align*}
which is a dense subspace of $\Cosc(\R^n)$. In particular, a continuous valuation on $\Cosc(\R^n)$ is uniquely determined by its restriction to $\Coco$.
Now, notice that any convex body in $\mathcal{K}^{n+1}$ induces a function in $\Coco$ by considering the lower part of its boundary as the graph of a convex function, and any function in $\Coco$ can be obtained in this way. More precisely, we associate to $K\in\K^{n+1}$ the function 
\begin{align*}
\lfloor K\rfloor(x)=\inf_{(x,t)\in K}t,
\end{align*}
which belongs to $\Coco$, compare Lemma \ref{lemma:lowerBoundaryInCoco}. In other words, we identify elements of $\Coco$ with convex bodies that coincide with the "lower part" of their epi-graphs. Consequently, any valuation on $\Coco$ induces a functional on $\K^{n+1}$ by precomposition with $\lfloor\cdot\rfloor$. In Section \ref{section:InducedValuations}, we examine the properties of this map and show that it is compatible with the different notions of translation invariance and the valuation property on the spaces $\Coco$ and $\K^{n+1}$. In addition, we show that it is continuous if we equip $\K^{n+1}$ with the Hausdorff metric. This allows us to obtain Theorem \ref{theorem:CharMaxDegree} from the following result by McMullen. Recall that a map $Y:\mathcal{K}^n\rightarrow\R$ is called $i$-homogeneous if
\begin{align*}
	Y(tK)=t^iY(K)\quad\text{for all }K\in\mathcal{K}^n, t\ge 0,
\end{align*}
and $Y$ is translation invariant if $Y(K+x)=Y(K)$ for every $K \in \K^n$ and $x \in \R^n$.
\begin{theorem}[McMullen \cite{McMullenContinuoustranslationinvariant1980}]\label{McM}
A functional $Y:\K^n \to \R$ is a continuous, translation invariant valuation which is $(n-1)$-homogeneous, if and only if there exists a continuous function $\eta:\sph^{n-1}\to \R$ such that 
\begin{equation}\label{e5}
Y(K)=\int_{\sph^{n-1}}\eta(\nu)\,dS_{n-1}(K,\nu)
\end{equation}
for every $K \in \K^{n}$. The function $\eta$ is uniquely determined up to the addition of the restriction of a linear function.
\end{theorem} 
Here, $S_{n-1}(K,\cdot)$ denotes the surface area measure of $K\in\mathcal{K}^n$. Note that any epi-translation invariant, continuous valuation on $\Coco$ that is epi-homogeneous of degree $n$ induces an $n$-homogeneous valuation on $\K^{n+1}$ under the above identification and thus admits a representation as in Theorem \ref{McM}. Assuming that the function $\eta\in C(\sph^n)$ is supported on the lower half sphere $\sph^n_-\coloneqq \{N\in \sph^n: N\cdot e_{n+1}<0\}$, a suitable change of coordinates reduces the representation in Theorem \ref{McM} to the representation in Theorem \ref{theorem:CharMaxDegree}. However, the function $\eta$ is only unique up to the restriction of linear functions, but this problem can be overcome with a rather simple argument. Consequently, we obtain the following result.
\begin{theorem}
    \label{theorem:CharacMaxDegCoco}
    If $Z:\Coco\rightarrow\R$ is a continuous, epi-translation invariant valuation which is epi-homogeneous of degree $n$, then there exists a unique continuous function $\eta:\sph^n\rightarrow\R$ with support compactly contained in $\sph^n_-$ such that
    \begin{align*}
        Z(\lfloor K\rfloor)=\int_{\sph^n} \eta(\nu) dS_{n}(K,\nu).
    \end{align*}
    for all $K\in\K^n$. In particular,
    \begin{align}\label{equation:McMforFunc}
        Z(u)=\int_{\dom(u)} \eta\left(\frac{(\nabla u(x),-1)}{\sqrt{1+|\nabla u(x)|^2}}\right)\sqrt{1+|\nabla u(x)|^2}\,dx
    \end{align}
    for all $u\in\Coco$, where $y\mapsto \eta\left(\frac{(y,-1)}{\sqrt{1+|y|^2}}\right)\sqrt{1+|y|^2}$ has compact support.
\end{theorem}
As $\Coco\subset\Cosc(\R^n)$ is dense, this provides, in fact, a new proof of Theorem \ref{theorem:CharMaxDegree}. However, we want to stress that this actually establishes a classification of valuations on the dense subspace $\Coco$. This can be shown with the original approach in \cite{ColesantiEtAlhomogeneousdecompositiontheorem2020} as well. 
Let us also remark that this restriction on the support of the functionals is a general property of this class of valuations, compare \cite{Knoerrsupportduallyepi2021}.\\

Our second main result provides a characterization of epi-translation invariant, continuous valuations on $\Cosc(\R^n)$ that are epi-homogeneous of degree $1$. It is motivated by the following classification result by Goodey and Weil \cite{GoodeyWeilDistributionsvaluations1984}.
\begin{theorem}[Goodey and Weil]\label{gw}
A functional $Y: \K^n \to \R$ is a continuous, translation invariant valuation which is homogeneous of degree $1$, if and only if there are two sequences $(L_j)_j,(W_j)_j$ of convex bodies in $\mathcal{K}^n$ such that 
\begin{equation}\label{one-hom}
 Y(K)=\lim_{j \to \infty} V(K,L_j,\dots,L_j)-V(K,W_j,\dots,W_j) 
\end{equation}
holds uniformly on compact subsets of $\mathcal{K}^{n}$.
\end{theorem}

Here, $V:(\K^n)^n \to \R$ denotes the so called so-called \textit{mixed volume}, which is a multilinear functional on $\mathcal{K}^n$ with respect to Minkowski addition (compare {\cite[Theorem 5.1.7]{SchneiderConvexBodiesBrunn2013}}). Equation \eqref{one-hom} can equivalently be written (compare {\cite[Chapter 5.1]{SchneiderConvexBodiesBrunn2013}}) as
\begin{align*}
	 Y(K)=&\lim_{j \to \infty} V(K,L_j,\dots,L_j)-V(K,W_j,\dots,W_j) \\
	 =&\lim_{j \to \infty} \frac{1}{n}\int_{\sph ^{n-1}}h_K \,dS_n(L_j)- \frac{1}{n}\int_{\sph ^{n-1}}h_K \,dS_n(W_j).
\end{align*}
Recall that Theorem \ref{theorem:CharacMaxDegCoco} is based on the observation that integrals with respect to the surface area measure are related to integrals involving the gradient of the associated convex functions by a suitable change of coordinates. Similarly, the Fenchel-Legendre transform of an element of $\Cosc(\R^n)$ may be considered as the support function of its epi-graph. Following these interpretations, we prove the following functional version of the result by Goodey and Weil.
\begin{theorem}
	\label{theorem:onehom}
	Every continuous, epi-translation invariant valuation $Z:\Cosc(\R^n)\rightarrow\R$ that is epi-homogeneous of degree $1$ can be approximated by a sequence $(Z_j)_j$ of valuations on $\Cosc(\R^n)$ that converges uniformly on compact subsets of $\Cosc(\R^n)$ with the following properties:
	\begin{enumerate}
		\item $Z_j$ is a continuous, epi-translation invariant valuation for each $j\in\mathbb{N}$.
		\item For every $j\in\mathbb{N}$ there exist two functions $\ell_j,w_j\in\Coco$ such that 
			\begin{align*}
				Z_j(u)=\int_{\dom(\ell_j)} u^{*}(\nabla \ell_j(x))\,dx-\int_{\dom(w_j)} u^{*}(\nabla w_j(x)) \,dx 
			\end{align*}
			for all $u\in\Coco$.
	\end{enumerate}
\end{theorem}
	Note that this representation only holds on the dense subspace $\Coco$, as the integrals are not convergent for general elements of $\Cosc(\R^n)$. This is a technical artifact of this representation, as the contributions of the two integrals cancel outside of a compact subset. More precisely, we deduce Theorem \ref{theorem:onehom} from the following result.
\begin{theorem}
    \label{theorem:onehomDual}
    Every continuous, epi-translation invariant valuation $Z:\Cosc(\R^n)\rightarrow\R$ that is epi-homogeneous of degree $1$ can be approximated uniformly on compact subsets of $\Cosc(\R^n)$ by a sequence $(Z_j)_j$ of valuations on $\Cosc(\R^n)$ given by
    \begin{align*}
        Z_j(u)=\int_{\R^n}\phi_j(x) u^*(x) dx,
    \end{align*}
    such that the following holds:
    \begin{enumerate}
        \item There exists $R>0$ such that $\phi_j$ is supported on $B_R(0)$ for all $j\in\mathbb{N}$.
        \item $\int_{\R^n}\phi_j(x)l(x)dx=0$ for all affine maps $l:\R^n\rightarrow\R$ and all $j\in\mathbb{N}$.
    \end{enumerate}
\end{theorem}
This result is implicitly contained in \cite{KnoerrSmoothvaluationsconvex}. However, we provide a simpler proof and show how one can explicitly obtain the sequence from a given valuation. In order to obtain Theorem \ref{theorem:onehom} from Theorem \ref{theorem:onehomDual}, we once again identify the valuations $Z_j$ with a valuation on convex bodies in $\R^{n+1}$. This reduces Theorem \ref{theorem:onehom} to showing that the function on $\sph^n$ obtained from $\phi_j$ can be written as a difference of surface area measures, which is a simple consequence of Minkowski's existence theorem.\\
 
    This paper is structured as follows. In Section \ref{section:Preliminaries}, we recall some facts from convex geometry, in particular the notion of support measures for convex bodies. In Section \ref{section:FromConvexBodiesToFunctions}, we discuss the relation between $\Coco$ and $\mathcal{K}^{n+1}$ and we examine how the support measures of a convex body in $\mathcal{K}^{n+1}$ are related to Hessian measures of the induced convex function. These ideas are used in Section \ref{section:InducedValuations} to show how valuations on $\Coco$ may be interpreted as valuations on $\mathcal{K}^{n+1}$. In Section \ref{section:ProofMaxDegree}, we use this correspondence to prove Theorem \ref{theorem:CharacMaxDegCoco}. Finally, in Section \ref{section:Theorem1Hom} we use a mollification procedure in the dual setting to prove Theorem \ref{theorem:onehomDual} and obtain Theorem \ref{theorem:onehom}.

\paragraph{Acknowledgements.} The authors would like to thank Monika Ludwig and Fabian Mussnig for many insightful comments and discussions. The second named author thanks the research unit Convex and Discrete Geometry at TU Wien, where he was hosted during part of this collaboration. 

\section{Preliminaries}
\label{section:Preliminaries}

We work in the Euclidean space $\R^{n+1}$, with $n\geq 1$. The vectors $e_1,\dots, e_{n+1}$ denote its standard orthonormal basis, and we identify $H\coloneqq \mathrm{span}\{e_1,\dots,e_n\}\subset \R^{n+1}$ with $\R^n$.
Both spaces are endowed with the standard Euclidean norm $|\cdot |$ and the standard scalar product $x\cdot y$ for $x, y \in \R^n$. We denote the ball with radius $R>0$ centered at $x\in \R^n$ by $B_R(x)$ or $B^n_R(x)$ if we want to emphasize its dimension. The $n$-dimensional Hausdorff measure is denoted by $\haus^n$.
\subsection{Convex bodies}
We denote by $\K^n$ the space of convex bodies in $\R^n$, that is, the set of all non-empty, convex, and compact subsets of $\R^{n}$, equipped with the Hausdorff metric. We refer to the monograph by Schneider \cite{SchneiderConvexBodiesBrunn2013} for a comprehensive background on convex bodies and only collect the notions and results we need in the following sections.\\
First, we require the following result on the convergence of intersections. Recall that two convex bodies $K, L\in\K^n$ are said to be separable by a hyperplane if there exists $y\in\R^n\setminus\{0\}$ and $c\in\R$ such that $K\subset \{x\in\R^n: \langle x,y\rangle\le c\}$, $L\subset \{x\in\R^n: \langle x,y\rangle\ge c\}$.
\begin{theorem}[\cite{SchneiderConvexBodiesBrunn2013} Theorem 1.8.10]
    \label{theorem:ConvergenceIntersectionNonSeperableBodies}
    Let $K, L\in\K^n$ be convex bodies that cannot be separated by a hyperplane. If $K_j,L_j\in \K^n$ are convex bodies with $K_j\rightarrow K$, $L_j\rightarrow L$ for $j\rightarrow\infty$, then $K_j\cap L_j\ne\emptyset$ for almost all $j\in\mathbb{N}$  and $K_j\cap L_j\rightarrow K\cap L$ for $j\rightarrow\infty$.
\end{theorem}
For every $K \in \K^{n}$, its \textit{support function} is the convex function $h_K:\R^n\rightarrow \R$ given by 
\begin{equation}\label{eq:support_function}
    h_K(x)=\sup_{y \in K}x \cdot y\quad \text{for}~x \in \R^n.
\end{equation}
Note that $K$ is uniquely determined by its support function. \\

Given a convex subset $A \subset \R^n$, we denote by $I_A:\R^n \to (-\infty,+\infty]$ the convex \textit{indicator function} of $A$, that is,
\begin{equation*}
I_A(x)= \begin{cases}
	  0 \qquad \qquad x \in A,\\
	  +\infty \qquad \text{otherwise}.
	 \end{cases}
\end{equation*}
Note that if $K \in \K^n$, then $I_K \in \Cosc(\R^n)$ and $(I_K)^*=h_K$.\\


Consider the set $\Sigma=\R^n \times \mathbb{S}^{n-1}$. For a fixed $K \in \K^n$, a pair $(x,\xi) \in \Sigma$ is called a \textit{support element} of $K$ if $x \in \partial K$ and $\xi$ is an outer unit normal to $\partial K$ in $x$, that is, $h_K(\xi)=x \cdot \xi$. We define the \textit{normal bundle} of $K$ as the set \[\text{Nor}K=\{(x,\xi) \in \Sigma\,:\, (x,\xi) \text{ is a support element of }K \}. \] 
For $x\in \R^n$, let $p(K,x)$ denote the point in $K$ with minimal distance to $x$, that is, $p(K,\cdot)$ is the nearest point projection. Note that for every $x \in \R^n\setminus K$, the pair $\left(p(K,x),\frac{x-p(K,x)}{|x-p(K,x)|}\right)$ is a support element of $K$. Consider the set $K_t$ of all points in $\R^n$ with distance smaller than $t>0$ from $K$ and the map 
\begin{align*}
 f_t\,: \,K_t \setminus K &\to \qquad\Sigma \\
 x &\mapsto \left(p(K,x),\frac{x-p(K,x)}{|x-p(K,x)|}\right),
\end{align*}
which is continuous and measurable. Let $\mu_t(K,\cdot)$ denote the measure on $\Sigma$ obtained by considering the push-forward of the Lebesgue measure on $K_t\setminus K$ along this map. 

We summarize the main properties of this measure in the following statement (see {\cite[Theorem 4.2.1]{SchneiderConvexBodiesBrunn2013}}).
\begin{theorem}\label{supportmeasures}
For every convex body $K \in \K^n$ there exist finite positive Borel measures $\Theta_i(K,\cdot),0\leq i \leq n-1,$ on $\mathcal{B}(\Sigma)$ such that for every $\beta \in \mathcal{B}(\Sigma)$ and every $t >0$
\[\mu_t(K,\beta)=\sum_{i=0}^{n-1}t^{n-i}\binom{n}{i}\Theta_i(K,\beta). \]

The mapping $K \mapsto \Theta_i(K,\cdot)$ is weakly continuous and is a measure-valued valuation.
\end{theorem}

 The coefficients $\Theta_i(K,\cdot)$ are called the \emph{support measures} of $K$. As an immediate consequence of Theorem \ref{supportmeasures}, we obtain the following: If $f \in C(\Sigma)$, the functional \begin{equation}\label{contval}
K \mapsto \int_\Sigma f(x,\xi)\,d\Theta_i\,(K,(x,\xi)) 
\end{equation}
is a continuous valuation for every $0\leq i \leq n-1$.  Note that $\Theta_i(K,\cdot)$ is concentrated on the normal bundle of $K$, compare \cite[p. 213]{SchneiderConvexBodiesBrunn2013}, so we can consider this as an integral over the normal bundle of $K$ instead of $\Sigma$.

The marginals of $\Theta_i(K,\cdot)$ with respect to $\R^n$ and $\sph^{n-1}$ are called the \textit{curvature measures} $C_i(K,\cdot)$ and \textit{area measures} $S_i(K,\cdot)$. If $K$ is a convex body with smooth boundary,  {\cite[Lemma 4.2]{SchneiderConvexBodiesBrunn2013}} implies that \eqref{contval} may be written as \[\int_{\partial K} f(x,\nu_K(x))\,dC_i(K, x), \] where $\nu_K:\partial K\rightarrow \sph^{n-1}$ is the \textit{Gauss map} of $K$, that is, $\nu_K(x)$ is the unique outer unit normal to $K$ in $x\in\partial K$. Similarly, if $K$ is strictly convex (that is, if $\partial K$ contains no segments) we may equivalently write \eqref{contval} as
\begin{equation}\label{continuousvaluation}
\int_{\mathbb{S}^{n-1}} f(\tau_K(\xi),\xi)\,dS_i(K, \xi),
\end{equation}
where $\tau_K(\xi)\in\partial K,\xi \in \sph^{n-1}$, is the unique boundary point such that $\xi$ is an outer normal. See \cite[Chapter 4]{SchneiderConvexBodiesBrunn2013} for more details.

\subsection{Convex functions}

\label{section:ConvexFunctions}

The spaces of functions used in this paper are all subsets of the space \[\Conv(\R^n)\coloneqq \{f: \R^n \to (-\infty,+\infty]: f \text{ is convex, proper and l.s.c.}\}. \] For a comprehensive background on convex functions, we refer to the monograph by Rockafellar and Wets \cite{RockafellarWetsVariationalAnalysis1998}. The space $\Conv(\R^n)$ is equipped with a natural topology induced by epi-convergence, where a sequence $(u_j)_j \subset \Conv(\R^n)$ is called \textit{epi-convergent} to $u \in \Conv(\R^n)$ if for every $x \in \R^n$ the following conditions hold: 
\begin{enumerate}
\item For every sequence $(x_j)_j$ that converges to $x$, $u(x)\leq \liminf_{j\to \infty} u_j(x_j)$.
\item There exists a sequence $(x_j)_j$ converging to $x$ such that $u(x)=\lim_{j \to \infty} u_j(x_j)$.
\end{enumerate}
In fact, this topology is metrizable, compare \cite[Theorem 7.58]{RockafellarWetsVariationalAnalysis1998}. 
For the space $\Cosc(\R^n)$, epi-convergence may be equivalently characterized by Hausdorff convergence of the level sets \begin{align*}
	\{u\leq t\}\coloneqq \{x \in \R^n: u(x)\leq t\}\quad\text{for }u\in\Cosc(\R^n), t\in\R.
\end{align*}
We will say that $\{u_j\leq t\} \to \emptyset$ if there exists $j_0$ such that $\{u_j\leq t\}=\emptyset$ for all $j\geq j_0$. The following result is a consequence of \cite[Lemma 5]{ColesantiEtAlValuationsConvexFunctions2017} and \cite[Theorem 3.1]{BeerEtAlCharacterizationEpiConvergence1992}.
\begin{lemma}\label{l1}
A sequence of functions $(u_j)_j \subset \Cosc(\R^n)$ epi-converges to $u \in \Cosc(\R^n)$ if and only if $\{u_j \leq t\} \to \{u \leq t\}$ for every $t \in \R$ with $t \neq \min_{x \in \R^n} u(x)$.
\end{lemma}

Let $u\in\Conv(\R^n,\R)$ and consider for $t>0$ and $\beta \in \mathcal{B}(\R^n\times \R^n)$ the set \[P_t(u,\beta)\coloneqq \{ x+ty\,:\, (x,y) \in \beta, y \in \partial u(x)\}, \]
where $\partial u(x)$ denotes the \emph{subdifferential} $u$ in $x$; that is, $y\in\partial u(x)$ if and only if $u(z)\ge u(x)+z\cdot (z-x)$ for all $z\in\R^n$. In other words, we are considering the image of the intersection of $\beta$ with the graph of the subgradient \[ \Gamma_u\coloneqq \{(x,y) \in \R^n\times \R^n\,:\, y \in \partial u(x) \}\]
under the map $\R^n\times\R^n\rightarrow\R^n$, $(x,y)\mapsto x+ty$. As shown by Colesanti and Hug \cite{ColesantiHugSteiner2000} for $u\in\Conv(\R^n,\R)$ and by Colesanti, Ludwig, and Mussnig \cite{ColesantiEtAlHessianValuations2020} in full generality,  the $n$-dimensional Hausdorff measure of $P_t(u,\beta)$ admits a Steiner-type formula similar to the support measures (compare Theorem \ref{supportmeasures}).
\begin{theorem}[Colesanti, Ludwig, and Mussnig {\cite[Theorem 7.1]{ColesantiEtAlHessianValuations2020}}]\label{theorem:HessianMeasures}
For $u \in \Conv(\R^n)$ there are non-negative Borel measures $\Xi_i(u,\cdot)$ on $\R^n \times \R^n$, $0\leq i \leq n$, such that \[\mathcal{H}^n(P_t(u,\beta))=\sum_{i=0}^n \binom{n}{i}t^{i} \Xi_{n-i}(u,\beta)  \] for every $\beta \in \mathcal{B}(\R^n\times \R^n)$ and $t>0$. These measures are called  \emph{Hessian measures}.
\end{theorem}
By {\cite[Theorem 7.3]{ColesantiEtAlHessianValuations2020}}, the Hessian measures are weakly-continuous with respect to epi-convergence. Moreover, if $u \in \Conv(\R^n)\cap C^2(\R^n)$ and $\beta \in \mathcal{B}(\R^n)$, then 
\begin{equation}\label{regular_hessian}
\Xi_i(u,\beta \times \R^n)=\int_{\beta} [\det D^2 u(x)]_{n-i}\, dx. 
\end{equation}

\section{From convex bodies to convex functions}
\label{section:FromConvexBodiesToFunctions}
\subsection{The space $\Coco$.}

Recall that our first result concerns valuations on the space \[\Coco\coloneqq \{u\in \Cosc(\R^n): u \text{ has compact domain} \}. \]

These functions can be obtained from convex bodies in $\R^{n+1}$ using the following construction: To every $K \in \K^{n+1}$ we associate the function $\lfloor K \rfloor:\R^n\rightarrow[-\infty,+\infty]$ defined by $\lfloor K \rfloor (x)=\inf\{t\in\R: (x,t)\in K\}$. Note that this is equivalent to
\begin{align*}
\lfloor K \rfloor (x)=\begin{cases}
      \min\{t\in\R: (x,t)\in K\} & x\in\mathrm{pr}_H(K),\\
      +\infty & \text{else},
  \end{cases}
\end{align*}
where $\mathrm{pr}_H:\R^{n+1}\rightarrow H\cong \R^n$ denotes the orthogonal projection.

\begin{lemma}
    \label{lemma:lowerBoundaryInCoco}
	$\lfloor K\rfloor\in \Coco$ for all $K\in\mathcal{K}^{n+1}$.
\end{lemma}
\begin{proof}
	Note first that $\lfloor K\rfloor$ is bounded from below by 
	\begin{align*}
		\inf\{t\in \R:(x,t)\in K\text{ for some }x\in\R ^n\},
	\end{align*} which is finite due to the compactness of $K$, so $\lfloor K\rfloor(x)\in (-\infty,+\infty]$ for $x\in\R^n$. Moreover, 
	\begin{align*}
		\dom \lfloor K\rfloor =\mathrm{pr}_H(K).
	\end{align*}
	In particular, $\dom \lfloor K\rfloor$ is compact and non-empty, so $\lfloor K\rfloor$ is proper.\\
	
	Let us show that $\lfloor K\rfloor$ is lower semi-continuous. If $x\in \dom(\lfloor K\rfloor)$ and $(x_j)_j$ is a sequence in $\dom(\lfloor K\rfloor)$ converging to $x$, then $(x_j,\lfloor K\rfloor(x_j))\in K$ for all $j\in\mathbb{N}$. This sequence is bounded in $\R^{n+1}$, so $t\coloneqq \liminf_{j\rightarrow\infty} \lfloor K\rfloor(x_j)$ exists and is finite. Thus, $(x,t)$ is a limit point of the sequence $(x_j,\lfloor K\rfloor(x_j))$ and therefore belongs to $K$. In particular,
	\begin{align*}
		\lfloor K\rfloor(x)\le t=\liminf_{j\rightarrow\infty} \lfloor K\rfloor(x_j).
	\end{align*}
	On the other hand, $x\in \R^n\setminus \dom(\lfloor K\rfloor)$ implies that $\lfloor K\rfloor$ is equal to $+\infty$ on a neighborhood of $x$, as the domain is closed. Thus $\lfloor K\rfloor$ is lower semi-continuous.\\
	It is easy to see that $\lfloor K\rfloor$ is convex. Thus, $\lfloor K\rfloor\in \Coco$ for all $K\in\mathcal{K}^{n+1}$.\\
\end{proof}

\begin{corollary}
	$\Coco\subset\Cosc(\R^n)$ is dense. 
\end{corollary}
\begin{proof}
	For $u\in \Cosc(\R^n)$, set $u_j\coloneqq \lfloor \epi(u)\cap (B_j^{n}(0)\times[-j,j])\rfloor$. Then $u_j\in \Coco$ for all $j\in\mathbb{N}$ large enough. As $u$ has compact sublevel sets, given $t\in\R$ we have
	\begin{align*}
		\{u_j\le t \}=\{u\le t \}\quad \text{for all }j\in\mathbb{N} \text{ large enough}.
	\end{align*}
	Lemma \ref{l1} thus implies that $(u_j)_j$ converges to $u$ in $\Cosc(\R^n)$, which shows the claim.\\
\end{proof}

\begin{lemma}
	\label{lemma:continutiyLowerBoundaryFunction}
	The map $\lfloor \cdot\rfloor:\mathcal{K} ^{n+1}\rightarrow\Coco$ is continuous.
\end{lemma}
\begin{proof}
	Let $(K_j)_j\subset \K^{n+1}$ be a sequence converging to $K \in \K^{n+1}$. We will use Theorem \ref{theorem:ConvergenceIntersectionNonSeperableBodies} to show that the relevant sublevel sets of $\lfloor K_j\rfloor$ converge to the corresponding sublevel sets of $\lfloor K\rfloor$. In order to simplify the argument, it is useful to consider instead the convex bodies
	\begin{align*}
		\tilde{K}_j\coloneqq K_j+[0,e_{n+1}],
	\end{align*} 
    which converge to $\tilde{K}\coloneqq K+[0,e_{n+1}]$. We may thus choose $R>0$ such that $\tilde{K}_j,\tilde{K}\subset B_R^{n}(0)\times [-R,R]$ for all $j\in\mathbb{N}$.\\
	
	As $\lfloor \tilde{K}_j\rfloor=\lfloor K_j\rfloor$, $\lfloor \tilde{K}\rfloor=\lfloor K\rfloor$ for all $j\in\mathbb{N}$, we obtain
	\begin{align*}
		\{\lfloor K_j\rfloor\le t\}=
		\mathrm{pr}_H(\tilde{K}_j\cap (B_R^{n}(0)\times [-(R+1),t])).
	\end{align*}
	and similarly for the sublevel sets of $\lfloor K\rfloor$. Note that the sets $\tilde{K}$ and $(B_R^{n}(0)\times [-(R+1),t]))$ can not be separated by a hyperplane for $t>\min_{(x,s)\in \tilde{K}}s=\min_{x\in\R^n}\lfloor K\rfloor(x)$. For $t>\min_{x\in\R^n}\lfloor K\rfloor(x)$, Theorem \ref{theorem:ConvergenceIntersectionNonSeperableBodies} thus implies 
	$\tilde{K}_j\cap (B_R^{n}(0)\times [-(R+1),t])\ne \emptyset$ for almost all $j\in\mathbb{N}$ and 
	\begin{align*}
		\tilde{K}_j\cap (B_R^{n}(0)\times [-(R+1),t])\rightarrow \tilde{K}\cap (B_R^{n}(0)\times [-(R+1),t])
	\end{align*}
	for $j\rightarrow\infty$. Applying the orthogonal projection onto $H$ to both sides, we obtain for $t>\min_{x\in\R^n}\lfloor K\rfloor(x)$
	\begin{align*}
		\{\lfloor K_j\rfloor\le t\}\rightarrow \{\lfloor K\rfloor\le t\}.
	\end{align*}
	On the other hand, $t<\min_{x\in\R^n}\lfloor K\rfloor(x)$ implies that $\{\lfloor K\rfloor\le t\}=\emptyset$. Therefore $\{\lfloor K_j\rfloor\le t\}=\emptyset$ for almost all $j\in\mathbb{N}$, as we may otherwise find a sequence $x_{j_k}\in \R^n$ with 
	\begin{align*}
		(x_{j_k},\lfloor K_{j_k}\rfloor(x_{j_k})) \in \tilde{K}_{j_k}\cap (B_R^{n}(0)\times [-(R+1),t]),
	\end{align*}
	from which we can construct a limit point $(x,t_0)\in \tilde{K}\cap (B_R^{n}(0)\times [-(R+1),t])$.
	
	Lemma \ref{l1} thus implies that $\lfloor K_j \rfloor\to\lfloor K \rfloor$. As $\Coco\subset\Cosc(\R^n)$ is metrizable, this shows that $\lfloor\cdot\rfloor$ is continuous.
\end{proof}

\subsection{On the relation of support measures and Hessian measures}

In \cite{ColesantiHugSteiner2000} Colesanti and Hug obtained the Hessian measures from the support measures of the epi-graph of $u\in\Conv(\R^n,\R)$ by pushing these measures to $\R^n\times\R^n$. This section is devoted to making this relation precise for elements of $\Coco$. However, we will not consider the support measures of the epi-graph of these functions and instead work directly with  convex bodies in $\R^{n+1}$ that define the same element on $\Coco$ via the map $\lfloor\cdot\rfloor$.\\
It will be advantageous to fix a suitable convex body inducing a given function $u\in\Coco$.

We set $M_u\coloneqq \max_{x \in \dom(u)} u(x)$, which is finite since the domain of $u$ is compact and $u$ is convex, and define \[K^u\coloneqq \epi (u-M_u) \cap R_{H}(\epi (u-M_u))+M_u e_{n+1}, \] where $R_{H}$ is the reflection with respect to $H$. Obviously, this is a compact and convex set, so $K^u \in \K^{n+1}$. We thus obtain a map
\begin{equation}\label{eq:function_to_body}
	\begin{gathered}
		\Coco \to \K^{n+1} \\ \qquad \qquad  u \mapsto K^u.
	\end{gathered}
\end{equation}
	By construction $u=\lfloor K^u \rfloor$ for $u\in\Coco$. In particular, $\lfloor\cdot\rfloor:\mathcal{K}^{n+1}\rightarrow\Coco$ is surjective and continuous. It is thus natural to ask whether the map in \eqref{eq:function_to_body}, which is a right inverse of $\lfloor \cdot\rfloor$, is also continuous. This is not the case, which can be shown using an example that we are going to exploit in the proof of Theorem \ref{theorem:CharacMaxDegCoco} below (which is in turn taken from \cite{ColesantiEtAlhomogeneousdecompositiontheorem2020}).\\

Consider the lower half-sphere $\sph^n_{-}\coloneqq \{X \in \sph^n: X\cdot e_{n+1}<0\}$. We define the \textit{lower boundary} of $K$ by 
\begin{align*}
	\partial K_{-}\coloneqq \{X \in \partial K: \text{some unit normal to $K$ in $X$ belongs to } \sph^n_{-}\}.
\end{align*}
In general, $\partial K_-$ is neither open nor closed. We can parametrize its closure using the map
\begin{align*}
	f_u:\dom (u) &\to \R^{n+1}\\
	x&\mapsto (x,u(x)),
\end{align*} 
for $u=\lfloor K\rfloor$. Note that $\partial K_-$ differs from its closure by a set of zero $n$-dimensional Hausdorff measure. 
If $\gamma:\R^{n+1} \times \sph_-^n \rightarrow\R$ is bounded and Borel measurable, we can apply the area formula to obtain
\begin{equation}\label{change1}
\int_{\partial {K^u_-}} \gamma(X, \nu_{K^u}(X))\,d\haus^n(X)=\int_{\dom(u)} \gamma\left((x,u(x)),\frac{(\nabla u(x),-1)}{\sqrt{1+|\nabla u(x)|^2}}\right)\sqrt{1+|\nabla u(x)|^2} \,dx, 
\end{equation}
where $\sqrt{1+|\nabla u(x)|^2}$ is the approximate Jacobian of $f_u$. Here we used that $(\nabla u(x),-1)$ is an outer normal vector to $\epi(u)$ in $(x,u(x))$ whenever $u$ is differentiable in $x$ (which is the case $\haus^n$-almost everywhere on the interior of $\dom(u)$).\\

If $u=\lfloor K\rfloor$ for a $C^2_+$ body $K$, then we have the following local description of the curvature measure of its epigraph. We include the proof for completeness.
\begin{lemma}\label{equivalence}
Let $K \in \K^{n+1}$ be $C^2_+$ body, $u\coloneqq \lfloor K\rfloor$. For every Borel set $\alpha \subset {\partial K}_{-}$, \[C_i(K,\alpha)=\int_{\alpha}\left( \frac{1}{\sqrt{1+|\nabla u(\mathrm{pr}_H X)|^2}}\right)^{n-i}[D^2 u(\mathrm{pr}_H X)]_{n-i}\,d\haus^n(X) \] for every $1\leq i \leq n$.
\end{lemma}
\begin{proof}
As $K$ is a $C^2_+$ body, the curvature measures are given by integrals of the elementary symmetric functions of the eigenvalues of the \textit{second fundamental form} $\sff_K$ of $K$, which is the bilinear form obtained from the differential of the Gauss map (see \cite[Section 2.5]{SchneiderConvexBodiesBrunn2013}). If $u=\lfloor K \rfloor$ then, for $(u(x),x)=X\in \partial K_-$, \[\nu_K((x,u(x)))=\frac{(\nabla u(x),-1)}{\sqrt{1+|\nabla u(x)|^2}}.\] 
If $V\in T_X\partial K$, then $V=(v,\nabla u(x)\cdot v)$ for some $v\in\R^n$ and thus $d\nu_K(X):\nu_K(X)^\perp=T_X\partial K\rightarrow T_{\nu_K(X)}\sph^n=\nu_K(X)^\perp$ is given by
\begin{align*}
    d\nu_K(X)V=\frac{(D^2 u(x)v,0)}{\sqrt{1+|\nabla u(x)|^2}}-\frac{\nabla u(x)^T D^2 u(x) v}{1+|\nabla u(x)|^2}\frac{(\nabla u(x),-1)}{\sqrt{1+|\nabla u(x)|^2}}.
\end{align*}
For $V_1,V_2\in T_X\partial K$, $V_i=(v_i,\nabla u(x)\cdot v_i)$ for $v_i\in\R^n$, the second fundamental form is therefore given by
\begin{align*}
    \sff_K(X)[V_1,V_2]=\frac{v_1^T D^2 u(x)v_2}{\sqrt{1+|\nabla u(x)|^2}}.
\end{align*}
In particular, the eigenvalues of $\sff_K(X)$ coincide with the eigenvalues of $D^2u(x)$ up to a factor of $\frac{1}{\sqrt{1+|\nabla u(x)|^2}}$, which implies the claim.
\end{proof}

Consider the \textit{gnomonic projection} $g:\sph^n\rightarrow\R^n$,
\[N \mapsto \frac{N-(N\cdot e_{n+1})N}{N\cdot e_{n+1}}. \]
If we set $\Sigma_-\coloneqq \R^{n+1}\times \sph^n_-$, we obtain a well-defined map
\begin{align*}
    (\mathrm{pr}_H,g):\Sigma_-&\rightarrow\R^n\times\R^n\\
    (X,N)&\mapsto (\mathrm{pr}_H(X),g(N)).
\end{align*}

The next result gives an explicit representation of support measures in terms of Hessian measures. The inverse relation was investigated by Colesanti and Hug in \cite[Theorem 3.1]{ColesantiHugSteiner2000}. In fact, Theorem \ref{theorem:RelationSupportMeasuresHessianMeasures} can be related to \cite[Theorem 3.1]{ColesantiHugSteiner2000} by a suitable change of coordinates using that if $N \in \sph^n_-$ is an outer normal vector to the epigraph of $u$ in $X=(x,u(x))$, and $y \in \partial u(x)$, one has \[ -N\cdot e_{n+1}=\frac{1}{\sqrt{1+|y|^2}}.\]
\begin{theorem}
    \label{theorem:RelationSupportMeasuresHessianMeasures}
    Let $K\in\K^{n+1}$, $u\coloneqq \lfloor K\rfloor\in\Coco$ and $0\leq i \leq n$. For every Borel set $\beta \subset \Sigma_-$, 
    \begin{align*}
        \Theta_i(K,\beta)=\int_{(\mathrm{pr}_H,g)(\beta)} \left(\frac{1}{\sqrt{1+|y|^2}} \right)^{n-i-1}\, d\Xi_i(u,(x,y)).
    \end{align*}
\end{theorem}
\begin{proof}
    Assume first that $K$ is a $C^2_+$ body.
    Let $\eta:\Sigma_-\rightarrow\R$ be bounded and measurable. Then, by \cite[Lemma 4.2.2]{SchneiderConvexBodiesBrunn2013}, Lemma \ref{equivalence}, and the area formula
    \begin{align*}
        &\int_{\Sigma_-} \eta(X,N)\,d\Theta_i(K,(X,N))\\ =&\int_{{\partial K}_{-}} \eta(X,\nu_K(X))\,dC_i(K,X)\\ =&\int_{{\partial K}_{-}} \eta(X,\nu_K(X))\left( \frac{1}{\sqrt{1+|\nabla u(\mathrm{pr}_H X)|^2}}\right)^{n-i}[D^2 u(\mathrm{pr}_H X)]_{n-i}\, d\haus^n(X)\\=&\int_{\dom(u)} \eta \left(x,u(x),\left(\frac{(\nabla u(x),-1)}{\sqrt{1+|\nabla u(x)|^2}} \right) \right)\left(\frac{1}{\sqrt{1+|\nabla u(x)|^2}} \right)^{n-i-1} [D^2 u(X)]_{n-1}\,dx.
    \end{align*}
    Then, by \eqref{regular_hessian}, \[\int_{\Sigma} \eta(X,N)\,d\Theta_i(K,(X,N))=\int_{\R^n\times \R^n}\eta((x,u(x)),y)\left(\frac{1}{\sqrt{1+|y|^2}} \right)^{n-i-1}\, d\Xi_i(u,(x,y)), \] which concludes the proof in the smooth case. For the general case, assume that $(K_j)_j$ is a sequence of $C^2_+$ bodies converging to $K$. By Lemma \ref{lemma:continutiyLowerBoundaryFunction}, the sequence $u_j\coloneqq \lfloor K_j\rfloor\in \Coco$ converges to $u$. As the support measures are continuous with respect to weak convergence, we thus obtain for $\eta\in C_c(\Sigma_-)$ 
    \begin{align*}
        &\int_{\Sigma_-} \eta(X,N)\,d\Theta_i(K,(X,N))\\
        =&\lim\limits_{j\rightarrow\infty}\int_{\Sigma_-} \eta(X,N)\,d\Theta_i(K_j,(X,N))\\
        =&\lim\limits_{j\rightarrow\infty}\int_{\R^n\times \R^n}\eta((x,u_j(x)),y)\left(\frac{1}{\sqrt{1+|y|^2}} \right)^{n-i-1}\, d\Xi_i(u_j,(x,y))\\
        =&\int_{\R^n\times \R^n}\eta((x,u(x)),y)\left(\frac{1}{\sqrt{1+|y|^2}} \right)^{n-i-1}\, d\Xi_i(u,(x,y)),
    \end{align*} 
    where we used \cite[Theorem 10.1]{ColesantiEtAlHessianValuations2020} in the last step. The claim follows.
\end{proof}

Following the  notation of \cite{ColesantiEtAlHadwigertheoremconvex2020}, for $u \in \Coco$ fixed we have the measures \[\Phi_i(u,\alpha)=\Xi_{n-i}(u;\alpha \times \R^n) \] and \[\Psi_i(u,\alpha)=\Xi_i(u;\R^n\times \alpha),\] whenever $0\leq i \leq n$ and $\alpha$ is a Borel subset of $\R^n$. The following corollary is a straightforward consequence of Theorem \ref{theorem:RelationSupportMeasuresHessianMeasures}.
\begin{corollary}
    Let $u\in\Coco$ and $0\leq i \leq n$.
    \begin{enumerate}
        \item If $K^u$ has non-empty interior and $\alpha\subset\R^{n+1}$ is a Borel subset of the graph of $u$ with $\alpha\subset \pi_H^{-1}(\mathrm{int}~\dom(u))$, then \[ C_i(K^u,\alpha)=\int_{\mathrm{pr}_H (\alpha)} \left(\frac{1}{\sqrt{1+|\nabla u(x)|^2}} \right)^{n-i-1}\, d\Phi_{n-i}(u;x).\]
        \item If $\alpha$ is a Borel subset of $\sph^n_{-}$, then \[ S_i(K^u,\alpha)=\int_{g(\alpha)} \left(\frac{1}{\sqrt{1+|y|^2}} \right)^{n-i-1}\, d\Psi_i(u;y).\]
    \end{enumerate}
\end{corollary}

The special case $i=n$ will be crucial later.
\begin{corollary}
	\label{lemma:change3}
	For every $u\in\Coco$, $\eta:\sph_-^n\rightarrow\R$ Borel measurable and bounded
	\begin{align}
		\label{eq:change3}
		\int_{\sph_-^n} \eta(N)\,dS_n(K^u,N)=\int_{\dom(u)} \eta\left(\frac{(\nabla u(x),-1)}{\sqrt{1+|\nabla u(x)|^2}}\right)\sqrt{1+|\nabla u(x)|^2} \,dx.
	\end{align}
    In particular, $x\mapsto \eta\left(\frac{(\nabla u(x),-1)}{\sqrt{1+|\nabla u(x)|^2}}\right)\sqrt{1+|\nabla u(x)|^2}$ is integrable over $\dom(u)$.
\end{corollary}
\begin{proof}
    This follows from Theorem \ref{theorem:RelationSupportMeasuresHessianMeasures} noting that
    \begin{align*}
        \int_{\R^n\times \R^n}\eta\left( \frac{(y,-1)}{\sqrt{1+|y|^2}}\right)\sqrt{1+|y|^2}\, &d\Xi_n(u,(x,y))\\
        =&\int_{\dom(u)} \eta\left(\frac{(\nabla u(x),-1)}{\sqrt{1+|\nabla u(x)|^2}}\right)\sqrt{1+|\nabla u(x)|^2} \,dx,
    \end{align*}
    compare \cite[Section 10.4]{ColesantiEtAlHessianValuations2020}.
\end{proof}
Note that 
\begin{align}
    \label{esupp}
    \lfloor K\rfloor^*(x)=h_K(x,-1),\quad x\in\R^n.
\end{align}
As the support function $h_K$ is $1$-homogeneous, we obtain the following corollary, which in particular implies that the integrals in Theorem \ref{theorem:onehom} are well-defined. 
\begin{corollary}
	\label{corollary:integrabilityCompactDomain}
	If $u,v\in \Coco$, then $x\mapsto v^*(\nabla u(x))$ is integrable on $\dom(u)$.
\end{corollary}
\begin{proof}
    By Corollary 
	\ref{lemma:change3}, the function 
    \begin{align*}
        x\mapsto h_{K^v}\left(\frac{\nabla u(x),-1)}{\sqrt{1+|\nabla u(x)|^2}}\right)\sqrt{1+|\nabla u(x)|^2}=h_{K^v}(\nabla u(x),-1)
    \end{align*}
    is integrable over $\dom (u)$, so the claim follows from $h_{K^v}(\cdot,-1)=\lfloor K^v\rfloor^*=v^*$.
\end{proof}

We remark that Corollary \ref{lemma:change3} already appeared, in weaker forms, in {\cite[Lemma 4.4]{ColesantiEtAlHadwigertheoremconvex2022a}} and {\cite[Lemma 8.1.4]{KnoerrThesis}}. These results can be deduced directly from Corollary \ref{lemma:change3} using \eqref{esupp}.

\subsection{From valuations on convex functions to valuations on convex bodies}
\label{section:InducedValuations}

We will use the map $\lfloor \cdot\rfloor:\mathcal{K}^{n+1}\rightarrow \Coco$ to interpret valuations on $\Coco$ as valuations on convex bodies in $\R^{n+1}$. This is based on the following simple observation.
\begin{lemma}
	\label{lemma:ValuationPropertyLowerBoundary}
	If $K,L \in \K^{n+1}$ are such that $K \cup L \in \K^{n+1}$, then 
	\begin{align*}
		&\lfloor K \cap L \rfloor=\lfloor K \rfloor \vee \lfloor L \rfloor, & \lfloor K \cup L \rfloor=\lfloor K \rfloor \wedge \lfloor L \rfloor.
	\end{align*}
\end{lemma}
\begin{proof}
	By definition
	\begin{align*}
		\lfloor K \cap L \rfloor(x)=&\inf\{t\in\R: (x,t)\in K\cap L\}\\
		\ge& \inf\{t\in\R: (x,t)\in K\}\vee\inf\{t\in\R: (x,t)\in K\cap L\}=\lfloor K\rfloor(x)\vee \lfloor L\rfloor(x),\\
		\lfloor K \cup L \rfloor(x)=&\inf\{t\in\R: (x,t)\in K\cup L\}\\
		\le& \inf\{t\in\R: (x,t)\in K\}\wedge\inf\{t\in\R: (x,t)\in K\cap L\}=\lfloor K\rfloor(x)\wedge \lfloor L\rfloor(x).
	\end{align*}
	On the other hand,
	\begin{align*}
		\dom(\lfloor K \cap L \rfloor)=\dom(\lfloor K \rfloor) \cap \dom(\lfloor L \rfloor),\quad\quad
		\dom(\lfloor K \cup L \rfloor)=\dom(\lfloor K \rfloor) \cup \dom(\lfloor L \rfloor),
	\end{align*}
	as the domains are just the image of the corresponding convex bodies under the projection onto $H\cong \R^{n}$. In particular, both sides of each of the inequalities are finite if and only if one of the two sides is finite. We thus only have to consider points belonging to the corresponding domains. Assume that $\lfloor K \cap L \rfloor(x)<+\infty$. As $\lfloor K\rfloor(x)\vee\lfloor L \rfloor(x)\le \lfloor K \cap L \rfloor(x)<+\infty$, 
	\begin{align*}
		\{(x,t)\in\R^{n+1}: t\in [\lfloor K\rfloor(x),\lfloor K\cap L\rfloor(x)]\}&\subset K,\\  \{(x,t)\in\R^{n+1}: t\in [\lfloor L\rfloor(x),\lfloor K\cap L\rfloor(x)]\}&\subset L
	\end{align*}
	by convexity, as the points corresponding to the boundary points belong to these sets. Thus $(x,\lfloor K \rfloor(x)\vee \lfloor L \rfloor(x))\in K\cap L$, which implies
	\begin{align*}
		\lfloor K\cap L\rfloor(x)\le \lfloor K \rfloor(x)\vee \lfloor L \rfloor(x).
	\end{align*} 
	Now assume that $\lfloor K\cup L\rfloor(x)<+\infty$. Then $(x,\lfloor K\cup L\rfloor(x))\in K\cup L$. Without loss of generality, we may assume that $(x,\lfloor K\cup L\rfloor(x))\in K$, and therefore
	\begin{align*}
		\lfloor K\rfloor(x)\wedge \lfloor L\rfloor (x)\le \lfloor K\rfloor(x)\le \lfloor K\cup L\rfloor(x)
	\end{align*}
	by the definition of $\lfloor K\rfloor(x)$.
\end{proof}

\begin{theorem}\label{theorem:PropertiesInducedValuations}
For $Z:\Coco\to \R$ consider $Y:\mathcal{K}^{n+1}\rightarrow\R$ defined by 
\begin{align*}
Y(K)=Z(\lfloor K\rfloor).
\end{align*}
Then $Y$ has the following properties:
\begin{enumerate}
\item If $Z$ is a valuation, then so is $Y$.
\item If $Z$ is continuous, then $Y$ is continuous with respect to the Hausdorff metric.
\item If $Z$ is epi-translation invariant, then $Y$ is translation invariant, that is
\begin{align*}
	Y(K+X)=Y(K)\quad \text{for all }K\in\mathcal{K}^{n+1},X\in\R^{n+1}.
\end{align*}
\item If $Z$ is epi-homogeneous of degree $j$, then $Y$ is $j$-homogeneous, that is,
\begin{align*}
	Y(tK)=t^{j}Y(K)\quad \text{for all }K\in\mathcal{K}^{n+1}, t>0.
\end{align*}
\end{enumerate}
\end{theorem}
\begin{proof} 
\begin{enumerate}
	\item If $K,L\in\mathcal{K}^{n+1}$ satisfy $K\cup L\in \mathcal{K}^{n+1}$, then 
		\begin{align*}
			&\lfloor K \cap L \rfloor=\lfloor K \rfloor \vee \lfloor L \rfloor \text{ and } \lfloor K \cup L \rfloor=\lfloor K \rfloor \wedge \lfloor L \rfloor.
		\end{align*}
		by Lemma \ref{lemma:ValuationPropertyLowerBoundary}. Thus
		\begin{align*}
			Y(K\cup L)+Y(K\cap L)=&Z(\lfloor K \cup L \rfloor)+Z(\lfloor K \cap L\rfloor)\\
			=&Z(\lfloor K \rfloor \wedge \lfloor L \rfloor)+Z(\lfloor K \rfloor \vee \lfloor L \rfloor)\\
			=&Z(\lfloor K\rfloor)+Z(\lfloor L\rfloor)=Y(K)+Y(L).
		\end{align*}
	\item If $Z$ is continuous, then $Y=Z\circ \lfloor \cdot\rfloor$ is continuous due to the continuity of $\lfloor\cdot\rfloor$, compare Lemma \ref{lemma:continutiyLowerBoundaryFunction}.
	\item For $X=(v,c)\in\R^{n}\times\R$ and $K\in\mathcal{K}^{n+1}$, the definition of $\lfloor K\rfloor$ implies for $x\in\R^n$
	\begin{align*}
		\lfloor K+X\rfloor(x)=&\inf\{s\in\R: (x,s)\in K+X\}=\inf\{s\in\R: (x-v,s-c)\in K\}\\
		=&\inf\{s+c: s\in\R, (x-v,s)\in K\}\\
		=&\lfloor K\rfloor(x-v)+c.
	\end{align*}
	If $Z$ is epi-translation invariant, we obtain
	\begin{align*}
		Y(K+X)=&Z(\lfloor K+X\rfloor)=Z(\lfloor K\rfloor(\cdot-v)+c)=Z(\lfloor K\rfloor)=Y(K).
	\end{align*}
	Thus $Y$ is translation invariant.
	\item For $t>0$ we calculate for $x\in\R^n$
	\begin{align*}
		\lfloor tK\rfloor(x)=&\inf\{s\in\R: (x,s)\in tK\}=\inf\left\{s\in\R: \left(\frac{x}{t},\frac{s}{t}\right)\in K\right\}\\
		=&\inf\left\{ts:s\in\R, \left(\frac{x}{t},s\right)\in K\right\}\\
		=&t\lfloor K\rfloor\left(\frac{x}{t}\right).
	\end{align*}
	Thus $\lfloor tK\rfloor=t\sq \lfloor K\rfloor$, which implies \[Y(tK)=Z(\lfloor tK \rfloor)=Z(t\sq \lfloor K \rfloor)=t^j Z(\lfloor K \rfloor)=t^j Y(K)\]
	whenever $Z$ is epi-homogeneous of degree $j$.
\end{enumerate}
\end{proof}

Let us conclude this section with a simple observation concerning Theorem \ref{theorem:PropertiesInducedValuations}, aimed at giving a better intuition of the relation between valuations on convex bodies and functions. As shown by Colesanti, Ludwig, and Mussnig in \cite[Theorem 1.1]{ColesantiEtAlHessianValuations2020}, any $\zeta\in C(\R\times \R^n\times \R^n)$ that has compact support with respect to the second and third variables defines a continuous valuation on $\Conv(\R^n)$ by 
$$Z_{\zeta,i}(u)=\int_{\R^n \times \R^n} \zeta(u(x),x,y)\, d\Xi_i(u,(x,y)).$$
In light of Theorem \ref{theorem:PropertiesInducedValuations}  and equation \eqref{contval}, this may be seen as a consequence of the properties of the support measures and their relation with the Hessian measures in Theorem \ref{theorem:RelationSupportMeasuresHessianMeasures}.
\section{Proof of Theorem \ref{theorem:CharacMaxDegCoco}}
\label{section:ProofMaxDegree}

In this section, we obtain the representation formula \eqref{equation:McMforFunc} as a consequence of Theorem \ref{McM}. In particular, this provides a new proof of the representation in Theorem \ref{theorem:CharMaxDegree}. More precisely, Theorem \ref{theorem:CharacMaxDegCoco} shows that the same representation holds for continuous, epi-translation invariant valuations on $\Coco$ that are epi-homogeneous of degree $n$. As $\Coco\subset\Cosc(\R^n)$ is dense, this directly establishes the representation formula for the corresponding space of valuations on $\Cosc(\R^n)$ by continuity.
\begin{proof}[Proof of Theorem \ref{theorem:CharacMaxDegCoco}]
Let $Z:\Coco\rightarrow\R$ be a continuous, epi-translation invariant valuation that is epi-homogeneous of degree $n$, and consider $Y:\K^{n+1}\rightarrow\R$ given by  $Y(K)\coloneqq Z(\lfloor K \rfloor)$. Then $Y$ is a valuation on $\K^{n+1}$ which is continuous, translation invariant, and $n$-homogeneous by Theorem \ref{theorem:PropertiesInducedValuations}. By McMullen's Theorem \ref{McM} there exists $\eta'\in C(\sph^n)$ such that \[Y(K)=\int_{\sph^{n}} \eta'(N)\,dS_n(K,N) \] for every $K \in \K^{n+1}$. If we define $\tilde{\eta}(N)\coloneqq [\eta'(N)+\eta'(R_H N) ]/2$, then the valuation \[\tilde{Y}(K)\coloneqq \int_{\sph^n}\tilde{\eta}(N)\,dS_n(K,N)\] thus satisfies \[Z(u)=Y(K^u)=\tilde{Y}(K^u). \]
We will work with $\tilde{Y}$ and the function $\tilde{\eta}\in C(\sph^n)$.\\

For a convex body $K$ in $H$ and $\ell>0$, consider the cylinder $C(K,\ell)=K\times [0,\ell]\in\mathcal{K}^{n+1}$. Then by definition $I_K=\lfloor C(K,\ell)\rfloor$, so 
\begin{align*}
	Z(I_K)=&\tilde{Y}(C(K,\ell))=2\tilde{\eta}(-e_{n+1})V_n(K)+ \int_{\sph^n \cap H} \tilde{\eta}(N)\,dS_n(C(K,\ell),N)\\
	=&2\tilde{\eta}(-e_{n+1})V_n(K)+\ell \int_{\sph^{n-1}} \tilde{\eta}(\nu)\,dS_{n-1}(K,\nu),
\end{align*} 
where we identify $\sph^{n-1}$ and $\sph^{n}\cap H$.
As the left-hand side of this equation is independent of $\ell>0$, we infer that 
\begin{equation}\label{E2}
	\int_{\sph^{n-1}} \tilde{\eta}(\nu)\,dS_{n-1}(K,\nu)=0
\end{equation}
for every $K \in \K^{n}$. We may consider the left-hand side of \eqref{E2} as a valuation on convex bodies in $\mathcal{K}^n$ that is continuous, translation invariant, and $(n-1)$-homogeneous. As it vanishes identically, McMullen's Theorem \ref{McM} implies that $\tilde{\eta}|_{\sph^n\cap H}$ is the restriction of a linear function to $\sph^n\cap H$.\\
In particular, there exists a linear function function $l:\R^{n+1}\rightarrow\R$ such that $\tilde{\eta}+l\equiv 0$ on the equator $\sph^n\cap H$, and we set $\hat{\eta}=\tilde{\eta}+\frac{1}{2}[l+l\circ R_H]$. Then $\hat{\eta}$ vanishes on the equator $\sph^n\cap H$. Using Corollary \ref{lemma:change3} and the fact that linear functions belong to the kernel of the surface area measure, we obtain for $u\in \Coco$
\begin{align*}
	Z(u)=&\tilde{Y}(K^u)=\int_{\sph^n}  \hat{\eta}(N)\,dS_n(K^u,N)=2\int_{\sph^n_{-}}  \hat{\eta}(N)\,dS_n(K^u,N)\\
	=&\int_{\dom(u)} 2 \hat{\eta}\left(\frac{(\nabla u(x),-1)}{\sqrt{1+|\nabla u(x)|^2}}\right)\sqrt{1+|\nabla u(x)|^2}\,dx,
\end{align*} which shows the representation in equation \eqref{equation:McMforFunc}. Here we used that $\hat{\eta}$ vanishes on the equator $\sph^n\cap H$ and is symmetric with respect to $H$. If we set
\begin{align*}
    \eta(y):=\begin{cases}
        2\hat{\eta}(N) & N\in S^n_-,\\
        0 & \text{else},
    \end{cases}
\end{align*}
then $\eta\in C(S^{n})$ is supported on the closure on $S^n_-$ and 
\begin{align*}
    Z(u)=\int_{S^{n}}\eta(N)dS_n(K^u,N)\quad\text{for }u\in\Coco.
\end{align*}
It remains to see that the support of $\eta$ is contained in $S^n_-$, or equivalently that the support of $\hat{\eta}$ has empty intersection with $H\cap \sph^n$. We will use an argument given by Colesanti, Ludwig, and Mussnig in the proof of {\cite[Proposition 27]{ColesantiEtAlhomogeneousdecompositiontheorem2020}}. Consider the function $\zeta(y)=2\hat{\eta}\left( \frac{(y,-1)}{\sqrt{1+|y|^2}}\right)\sqrt{1+|y|^2}$ defined on $\R^n$. Then it is sufficient to show that $\zeta$ is compactly supported. Suppose by contradiction that the support is not compact. Then we can find a sequence $y_j\in\R^n$ such that $|y_j|\to \infty$, $\zeta(y_j)\neq 0$ for every $j\in\mathbb{N}$ and 
\[\lim_{j \to \infty} \frac{y_j}{|y_j|}=\nu\in\sph^n. \]
 Consider the sets 
 \[B_j\coloneqq \{x \in y_j^{\perp}: |x|\leq 1 \}, \quad B_{\infty}\coloneqq \{x \in \nu^\perp: |x|\leq 1 \} \] 
 and define the cylinders 
 \[C_j\coloneqq \left\{x+t\frac{y_j}{|y_j|}: x\in B_j, \ t \in \left[0, \frac{1}{|\zeta(y_j)|} \right] \right\}. \] 

For $y \in \R^n$ let $l_y$ denote the linear function $x\mapsto x\cdot y$. Consider the sequence \[u_j=l_{y_j}+I_{C_j} \] in $\Coco$. By construction, $\min_{x\in\dom(u_j)}u_j(x)=0$. For $t>0$, the sublevel sets are given by
\begin{align*}
	\{u_j\le t\}=\left\{x+t\frac{y_j}{|y_j|}: x\in B_j, \ s \in \left[0, \min\left\{\frac{t}{|y_j|},\frac{1}{|\zeta(y_j)|}\right\} \right] \right\},
\end{align*}
so $\{u_j\le t\}\rightarrow B_\infty$ in this case. Obviously, the sublevel sets are empty for $t<0$. Lemma \ref{l1} thus implies that $(u_j)_j$ converges to $I_{B_\infty}$.\\

Now note that $S_n(K^{u_j},\cdot)$ is concentrated on $(H\cap \sph^n)\cup \left\{\frac{(y_j,-1)}{\sqrt{1+|y_j|^2}},\frac{(y_j,1)}{\sqrt{1+|y_j|^2}}\right\}$, so
 \begin{align*}
 	Z(u_j)=&\int_{\sph^n}\hat{\eta}(N)\,dS_n(K^{u_j},N)=\left[\hat{\eta}\left(\frac{(y_j,-1)}{\sqrt{1+|y_j|^2}}\right)+\hat{\eta}\left(\frac{(y_j,1)}{\sqrt{1+|y_j|^2}}\right)\right] \sqrt{1+|y_j|^2}\vol_n(C_j)\\
 	=&\zeta(y_j)\vol_n(C_j)=\kappa_{n-1}, 
 \end{align*} 
because $\hat{\eta}$ is symmetric with respect to $H$ and vanishes on $\sph^n\cap H$. Here, $\vol_n$ denotes the $n$-dimensional Lebesgue measure and $\kappa_{n-1}$ is the volume of the $(n-1)$-dimensional unit ball. By continuity we obtain \[Z(I_{B_\infty})=\lim_{j\to \infty}Z(u_j)=\kappa_{n-1}.\]
On the other hand, $Z$ is $n$-homogeneous and $B_\infty$ is a convex body of dimension $n-1$, so $Z(I_{B_\infty})=0$, which is a contradiction. Thus $\zeta$ has compact support. \\

Finally, let us show how one can use McMullen's Theorem \ref{McM} to see that $\eta$ is uniquely determined by the valuation $Z$. Let us thus assume that $\eta,\eta'\in C(\sph^n)$ are compactly supported on $\sph^n_-$ and satisfy
\begin{align*}
	Z(\lfloor K\rfloor )&=\int_{\sph^n}\eta(N) \,dS_n(K,N)=\int_{\sph^n}\eta'(K) \,dS_n(K,N)\quad\text{for}~K\in\K^{n+1}.
\end{align*}
By McMullen's Theorem \ref{McM}, $\eta$ and $\eta'$ thus differ by the restriction of a linear function to $\sph^n$. However, they are both equal to $0$ on the complement of $\sph_-^n$, and therefore $\eta-\eta'$ vanishes on an open subset. As this difference is the restriction of a linear function, it thus has to vanish identically. In particular, $\eta=\eta'$.
\end{proof}

Let us add the following observation.
\begin{corollary}
	Let $Z:\Coco\rightarrow\R$ be a continuous and epi-translation invariant valuation that is epi-homogeneous of degree $n$. Then $Z$ extends uniquely to a continuous valuation on $\Cosc(\R^n)$.
\end{corollary}
\begin{proof}
	By Theorem \ref{theorem:CharacMaxDegCoco}, any such valuation $Z:\Coco\rightarrow\R$ is given by
	\begin{align*}
		Z(u)=\int_{\dom(u)}\zeta(\nabla u(x))\,dx\quad\text{for }u\in\Coco
	\end{align*}
	for some $\zeta\in C_c(\R^n)$. The right-hand side of this equation defines a continuous valuation on $\Cosc(\R^n)$ by \cite[Theorem~1.2]{ColesantiEtAlhomogeneousdecompositiontheorem2020}, which yields the desired continuous extension. As $\Coco \subset \Cosc(\R^n)$ is dense, this extension is unique.
\end{proof}
It is an interesting question whether this result holds for arbitrary continuous and epi-translation invariant valuations on $\Coco$, and we will address this problem in a future work.
\section{Proof of Theorem \ref{theorem:onehom}}
\label{section:Theorem1Hom}
For the proof of Theorem \ref{theorem:onehom}, we will switch to the dual setting: Recall that for any functional $Z:\Cosc(\R^n)\rightarrow\R$, we may define a functional $Z^*:\Conv(\R^n,\R)\rightarrow\R$ by
\begin{align*}
	Z^*(u)\coloneqq Z(u^*)\quad\text{for }u\in\Conv(\R^n,\R),
\end{align*}
where $u^*$ denotes the Fenchel-Legendre transform. Then the following holds, compare the discussion in {\cite[Section 3.1]{ColesantiEtAlHessianValuations2020}}.
\begin{itemize}
	\item $Z$ is a valuation if and only if $Z^*$ is a valuation.
	\item $Z$ is continuous if and only if $Z^*$ is continuous.
	\item $Z$ is epi-translation invariant if and only if $Z^*$ is \textit{dually epi-translation invariant}, that is, invariant with respect to the addition of affine functions to its argument.
	\item $Z$ is epi-homogeneous of degree $i$ if and only if $Z^*$ is $i$-homogeneous in the classical sense, that is,
	\begin{align*}
	Z^*(tu)=t^iZ^*(u)\quad\text{for all }u\in\Conv(\R^n,\R), t\ge 0.
	\end{align*}
\end{itemize}
Now assume that $Z$ is epi-homogeneous of degree $1$. According to {\cite[Corollary 24]{ColesantiEtAlhomogeneousdecompositiontheorem2020}}, $Z^*$ is then an \emph{additive} valuation, that is,
\begin{align*}
	Z^*(u+v)=Z^*(u)+Z^*(v)\quad\text{for all }u,v\in\Conv(\R^n,\R).
\end{align*}
In \cite{Knoerrsupportduallyepi2021}, this property was used to lift dually epi-translation invariant valuations to \textit{distributions}, that is, continuous linear functionals on the space of smooth functions with compact support. We refer to \cite{HoermanderLPDO} for a background on distributions.
\begin{theorem}[\cite{Knoerrsupportduallyepi2021} Theorem 2]
	For every $1$-homogeneous, dually epi-translation invariant, continuous valuation $Z:\Conv(\R^n,\R)\rightarrow\R$ there exists a unique distribution $\GW(Z)$ on $\R^n$ with compact support which satisfies
	\begin{align}
	\label{eq:PropertyGWDistribution}
	\GW(Z)[u]=Z(u) \quad\text{for all }u\in\Conv(\R^n,\R)\cap C^\infty(\R^n).
	\end{align}
\end{theorem}
Note that \eqref{eq:PropertyGWDistribution} is well-defined due to the compactness of the support. Let us remark that a similar result holds for homogeneous valuations of arbitrary degree of homogeneity, which is based on ideas of Goodey and Weil \cite{GoodeyWeilDistributionsvaluations1984} for translation invariant valuations on convex bodies. We refer to \cite{Knoerrsupportduallyepi2021} for this more general construction.

For a $1$-homogeneous, dually epi-translation invariant, and continuous valuation $Z$ on $\Conv(\R^n,\R)$, we define its \textit{support} as \[\supp Z\coloneqq \supp \GW(Z).\] This is a compact subset of $\R^n$ which has the property that \begin{align*}
	Z(u)=Z(v)\quad\text{for all }u,v\in\Conv(\R^n,\R)\text{ s.t. }u\equiv v\text{ on a neighborhood of }\supp Z,
\end{align*}
compare {\cite[ Proposition 6.3]{Knoerrsupportduallyepi2021}}. In particular, the support of such a valuation is a compact subset of $\R^n$.\\

We are now able to prove the following more general version of Theorem \ref{theorem:onehomDual}. Let us remark that Theorem \ref{theorem:onehomDual} can also be deduced from {\cite[Theorem 1]{KnoerrSmoothvaluationsconvex}} in combination with Lemma 5.3 of the same article. The proof we give here is self-contained and does not rely on the machinery developed in \cite{KnoerrSmoothvaluationsconvex}. In addition, it provides an explicit way to construct the approximating sequence.
\begin{theorem}\label{theorem:approx1hom}
		Let $Z: \Conv(\R^n,\R) \to \R$ be a valuation that is continuous, dually epi-translation invariant, homogeneous of degree 1, and satisfies $\supp\mu\subset A$ for some compact convex subset $A\subset\R^n$. Let $\phi\in C^\infty_c(\R^n)$ be a non-negative function with $\int_{\R^n}\phi(x)\,dx=1$ and $\supp\phi\subset B_1(0)$, and set
        \begin{align*}
            \phi_j(x)\coloneqq j^{n}\GW(\mu)\left[\phi\left(j(\cdot-x)\right)\right]\quad\text{for}~j\in\mathbb{N}. 
        \end{align*}
        Then the following holds:
        \begin{enumerate}
            \item $\phi_j\in C^\infty_c(\R^n)$,
			\item $\supp\phi_j\subset A+\frac{1}{j}B_{1}(0)$ for all $j\in\mathbb{N}$,
			\item $\int_{\R^n}\phi_j(x)l(x)\,dx=0$ for all affine functions $l:\R^n\rightarrow \R$ and all $j\in\mathbb{N}$.
		\end{enumerate}
        Moreover, the continuous, dually epi-translation invariant valuations $Z_j$ given by
	\begin{align*}
		Z_j(v)\coloneqq \int_{\R^n}v(x)\phi_j(x)\,dx
	\end{align*}
	converge uniformly to $Z$ on compact subsets of $\Conv(\R^n,\R)$.
	\end{theorem}
	\begin{proof}
		Let $T\coloneqq \GW(Z)$ denote the Goodey-Weil distribution of $Z$. Fix a non-negative function $\phi\in C^\infty_c(\R^n)$ with $\int_{\R^n}\phi(x)\,dx=1$, $\supp\phi\subset B_1(0)$ and consider the convolution $T_j\coloneqq T*j^{n}\phi(j\cdot)$ defined by
		\begin{align*}
			T_j(\psi)=T\left(\psi*j^{n}\phi\left(j\cdot\right)\right)=T\left(\int_{\R^n}\psi(x)j^{n}\phi\left(j(\cdot-x)\right)\,dx\right)\quad \text{for } \psi\in C^\infty_c(\R^n).
		\end{align*}
		As $T$ is continuous, we thus obtain 
		\begin{align*}
			T_j(\psi)=\int_{\R^n}\psi(x)j^{n}T\left(\phi\left(j(\cdot-x)\right)\right)\,dx\quad \text{for } \psi\in C^\infty_c(\R^n).
		\end{align*}
		Note that by definition,
		\begin{align*}
			\phi_j(x)=j^{n}T\left(\phi\left(j(\cdot-x)\right)\right).
		\end{align*}
		Elementary facts about the convolution of distributions (compare \cite[Chapter 4]{HoermanderLPDO}) show that $\phi_j\in C^\infty_c(\R^n)$ with $\supp\phi_j\subset A+\frac{1}{j}B_1(0)$ for all $j\in\mathbb{N}$. If $l\in C^\infty(\R^n)$ is an affine function, then
		\begin{align*}
			\left[l*j^{n}\phi\left(j\cdot\right)\right](y)=\int_{\R^n}l(y-x)j^{n}\phi\left(jx\right)\,dx
		\end{align*}
		is affine as well. As $T$ has compact support and $Z$ vanishes on affine functions, we obtain
		\begin{align*}
			\int_{\R^n}l(x)\phi_j(x)\,dx=T\left(l*j^{n}\phi\left(j\cdot\right)\right)=Z\left(l*j^{n}\phi\left(j\cdot\right)\right)=0,
		\end{align*}
		which shows the second property. In particular, if we define $Z_j:\Conv(\R^n,\R)\rightarrow\R$ by \begin{align*}
			Z_j(v)\coloneqq \int_{\R^n}v(x)\phi_j(x)\,dx,
		\end{align*}
		then $Z_j$ is a continuous, dually epi-translation invariant valuation that satisfies $\GW(Z_j)=T_j$. It remains to check that $(Z_j)_j$ converges to $Z$ uniformly on compact subsets. To see this, note that for $v\in\Conv(\R^n,\R)\cap C^\infty(\R^n)$ the function
		 \begin{align*}
		 	\left[v*j^{n}\phi\left(j\cdot\right)\right](y)=\int_{\R^n}v(y-x)j^{n}\phi\left(jx\right)\,dx,
		 \end{align*}
		 is convex as $\phi$ is non-negative. Thus, the fact that $T$ has compact support implies that 
		 \begin{align*}
		 	T\left(v*j^{n}\phi\left(j\cdot\right)\right)=\int_{\R^n}T(v(\cdot-x))j^{n}\phi\left(jx\right)\,dx.
		 \end{align*}
		 In particular, for any $v\in\Conv(\R^n,\R)\cap C^\infty(\R^n)$
		 \begin{align*}
		 	Z_j(v)=&\GW(Z_j)[v]=T_j(v)=\int_{\R^n}T(v(\cdot-x))j^{n}\phi\left(jx\right)\,dx\\
		 	=&\int_{\R^n}Z(v(\cdot-x))j^{n}\phi\left(jx\right)\,dx.
		 \end{align*}
	 	On $\Conv(\R^n,\R)$, the topology induced by epi-convergence coincides with the topology of uniform convergence on compact subsets of $\R^n$, see {\cite[Theorem 7.17]{RockafellarWetsVariationalAnalysis1998}}. Using this fact, it is easy to see that the map
	 	\begin{align*}
	 		\Conv(\R^n,\R)\times \R^n&\rightarrow\Conv(\R^n,\R)\\
	 		(v,x)&\mapsto v(\cdot -x)
	 	\end{align*}
	 	is continuous. In particular, 
		 \begin{align*}
		 	v\mapsto \int_{\R^n}Z(v(\cdot-x))j^{n}\phi\left(jx\right)\,dx
		 \end{align*} defines a continuous valuation on $\Conv(\R^n,\R)$. By continuity, we thus obtain
		 \begin{align*}
		 	Z_j(v)=\int_{\R^n}Z(v(\cdot-x))j^{n}\phi\left(jx\right)\,dx\quad\text{for all }v\in\Conv(\R^n,\R).
		 \end{align*}
		 Let $\epsilon>0 $ be given. Our previous discussion implies
		 \begin{align*}
		 	\left|Z_j(v)-Z(v)\right|=&\left|\int_{\R^n}Z(v(\cdot-x))j^{n}\phi\left(jx\right)\,dx-\int_{\R^n}Z(v)j^{n}\phi\left(jx\right)\,dx\right|\\
		 	\le &\int_{\R^n} \left|Z(v(\cdot-x))-Z(v)\right|j^{n}\phi\left(jx\right)\,dx.
		 \end{align*} 
	 	As the map
	 	\begin{align*}
	 		\Conv(\R^n,\R)\times \R^n&\rightarrow\Conv(\R^n,\R)\\
	 		(v,x)&\mapsto v(\cdot -x)
	 	\end{align*}
 		is continuous, it is uniformly continuous on compact subsets. Given a compact subset $K\subset \Conv(\R^n)$, we can thus find $\delta>0$ such that 
		 \begin{align*}
		 	|Z(v(\cdot-x))-Z(v)|\le \epsilon \quad\text{for all }v\in K\text{ and all } x\in\R^n\text{ with }|x|<\delta.
		 \end{align*}
		 As $\phi$ is supported on $B_1(0)$, $\supp \phi\left(j\cdot\right)\subset B_\delta(0)$ for all $j\ge \frac{1}{\delta}$, so
		 \begin{align*}
		 	\left|Z_j(v)-Z(v)\right|\le\int_{\R^n} \left|Z(v(\cdot-x))-Z(v)\right|j^{n}\phi\left(jx\right)\,dx\le \epsilon\quad \text{for all }v\in K \text{ and } j\ge \frac{1}{\delta}.
		 \end{align*}
		 Thus $(Z_j)_j$ converges uniformly to $Z$ on the compact subset $K\subset \Conv(\R^n,\R)$, which concludes the proof.
	\end{proof}

\begin{proof}[Proof of Theorem \ref{theorem:onehom}]
	Let $Z:\Cosc(\R^n)\rightarrow\R$ be a continuous, epi-translation invariant valuation that is epi-homogeneous of degree $1$ and assume that $\supp Z^*\subset B_R(0)$.\\
	By Theorem~\ref{theorem:approx1hom}, there exists a sequence $\phi_j\in C^\infty_c(\R^n)$ with $\supp \phi_j\subset B_{R+1}(0)$ such that 
	\begin{align*}
		Z_j(u)\coloneqq \int_{\R^n}u^*(x)\phi_j(x)\,dx
	\end{align*}
	defines a sequence of continuous, epi-translation invariant, and epi-homogeneous valuations of degree $1$ that converges uniformly to $Z$ on compact subsets. Here we use that the Fenchel-Legendre transform establishes a homeomorphism between $\Cosc(\R^n)$ and $\Conv(\R^n,\R)$, compare \cite[Theorem~11.34]{RockafellarWetsVariationalAnalysis1998}, so the preimage of any compact subset of $\Conv(\R^n,\R)$ under this map is compact.\\
	
	It remains to see that $Z_j$ has the desired representation on the subspace $\Coco$. Consider the function $b=\lfloor B^{n+1}_1(0)\rfloor \in \Coco$, that is,
	\begin{align*}
	b(x)=\begin{cases}
	1-\sqrt{1-|x|^2}& |x|\le 1,\\
	+\infty &|x|>1.
	\end{cases}
	\end{align*}
	Then $b^*(x)=\sqrt{1+|x|^2}$. From a direct calculation one infers that $\det D^2 b^{*}(x)=(1+|x|^2)^{-(n/2 + 1)}$, and using \eqref{change1} we can thus write for $u\in\Coco$ 
	\begin{align*}
	Z_j(u)=&\int_{\R^n}u^{*}(x)\phi_j(x)\,dx=\int_{\R^n}u^{*}(x)\phi_j(x)(1+|x|^2)^{n/2 + 1}\det D^2b^{*}(x)\,dx\\
	=&\int_{\dom(b)} u^{*}(\nabla b(x))\phi_j(\nabla b(x))(1+|\nabla b(x)|^2)^{n/2+1}\,dx=\int_{\sph^{n}_{-}} \frac{u^{*}(g(N))}{\sqrt{1+|g(N)|^2}}f_j(N)dN,
	\end{align*} 
	where \begin{align*}
			g: \sph^n_{-} &\to H\cong\R^n \\
			N &\mapsto \frac{N}{N\cdot e_{n+1}}-e_{n+1}
	\end{align*}
	and $f_j(N)\coloneqq \phi_j(g(N))(1+|g(N)|^2)^{n/2+2}$ is a function that is compactly supported on the lower half sphere.	We trivially extend $f_j$ to a smooth function on $\sph^{n}$.\\
	
	By equation \eqref{esupp}, we thus obtain the representation \[ Z_j(u)=\int_{\sph^n} h_{K^u}(N)f_j(N)dN.\]
	In fact, $Z_j(u)=\int_{\sph^n} h_{K}(N)f_j(N)dN$ for any $K\in\mathcal{K}^{n+1}$ with $h_K=h_{K^u}$ on $\sph^n_-$. As $Z_j$ is epi-translation invariant and epi-homogeneous of degree $1$, we thus obtain
	\begin{align*}
	0=Z_j( I_{\{v\}}+c)=\int_{\sph^n} h_{\{(v,c)\}}(N)f_j(N)dN\quad \text{for all }(v,c)\in \R^{n+1}.
	\end{align*}
	As $h_{\{(v,c)\}}(N)=(v,c)^T\cdot N$, the non-negative measure 
	\begin{align*}
	\mu_j(B)\coloneqq \int_B (1+\|f_j\|_\infty+f_j)dN \quad\text{for a Borel subset }B\subset\sph^n
	\end{align*}
	is thus not concentrated on a great sphere and satisfies 
	\begin{align*}
	\int_{\sph^n}N \,d\mu_j(N)=0.
	\end{align*} 
	By Minkowski's existence Theorem (see {\cite[Section 8.2.1]{SchneiderConvexBodiesBrunn2013}}), there thus exists a convex body $L_j\in \mathcal{K}^{n+1}$ such that $\mu_j=S_{n}(L_j)$. In particular,
	\begin{align*}
	Z_j(u)=\int_{\sph^n}h_{K^u} f_j(N)dN=\int_{\sph^n}h_{K^u} \,dS_{n}(L_j)-\int_{\sph^n}h_{K^u} \,dS_{n}(\sqrt[n]{1+\|f_j\|_\infty}B_1(0)).
	\end{align*}
	Here we have used that the surface area measure on $\R^{n+1}$ is $n$-homogeneous and that $S_n(B_1(0))$ is the spherical Lebesgue measure.	Set $W_j\coloneqq \sqrt[n]{1+\|f_j\|_\infty}B_1(0)$. By construction, $S_n(L_j)$ and $S_n(W_j)$ are absolutely continuous with respect to the spherical Lebesgue measure, and their densities only differ on the support of $f_j$, that is, on a compact subset contained in the lower half sphere. Set $\ell_j=\lfloor L_j \rfloor, w_j=\lfloor W_j\rfloor$. Applying Corollary \ref{lemma:change3} again and using that support functions are $1$-homogeneous, we thus obtain for $u\in \Coco$
	\begin{align*}
	Z_j(u)=&\int_{\sph^n_-}h_{K^u} \,dS_{n}(L_j)-\int_{\sph_-^n}h_{K^u} \,dS_{n}(W_j)\\
	=&\int_{\dom(\ell_j)} u^{*}(\nabla \ell_j(x)) \,dx-\int_{\dom(w_j)}u^{*}(\nabla w_j(x))\,dx.
	\end{align*}
	Thus $Z_j$ has the desired representation. In particular, $Z$ can be approximated uniformly on compact subsets by valuations of this type.
\end{proof}

Using the main results of \cite{Knoerrsupportduallyepi2021}, it is possible to deduce Theorem \ref{theorem:onehom} from the characterization of $1$-homogeneous valuations on convex bodies by Goodey and Weil in Theorem \ref{one-hom} using the ideas in Section \ref{section:FromConvexBodiesToFunctions}. To avoid unnecessary technicalities, we have decided to present a more direct argument.

\bibliography{biblio2} 
\bibliographystyle{siam}

\begin{tabular}{l l}
Jonas Knoerr & Jacopo Ulivelli\\
Institute of Discrete Mathematics and Geometry & Diparitimento di Matematica\\
TU Wien & Sapienza, University of Rome\\
Wiedner Hauptstrasse 8-10, 1040 Wien & Piazzale Aldo Moro 5, 00185 Rome\\
Austria & Italy \\
  jonas.knoerr@tuwien.ac.at & jacopo.ulivelli@uniroma1.it\\
\end{tabular}

\end{document}